\documentclass[12pt]
{article}
\usepackage{amssymb,amsmath}

\usepackage[latin1]{inputenc}
\usepackage[T1]{fontenc}
\usepackage{play}

\usepackage{graphicx}
\DeclareGraphicsExtensions{.pdf,.png,.jpg}

\graphicspath{ {images/} }

\begin{document}
\title{
Guess the Larger Number
} 
\author {Alexander Gnedin\footnote{\tt a.gnedin@qmul.ac.uk}} 
\maketitle

{\small 
\begin{play}{
\speaker{B} 
 I myself have invented a game. Well, think of a number.
\speaker{A} I got a number.
\speaker{B} Me  too. Now, tell me yours.
\speaker{A} Seven.
\speaker{B} Seven. Mine is eight -- I won.}
\vskip0.2cm
\speaker {\rm Sergey Solovyov,} {\it Assa} (conversation of Bananan and Alika)  
 \end{play}
}

\begin{abstract}
\noindent
We discuss variations of the zero-sum game where Bob  
selects two distinct numbers,  and Alice learns one of them to make a guess 
which of the numbers is the larger. 
\end{abstract}

\noindent

\newpage

\noindent
\section{Introduction.}
To play a guessing  game,  Bob   takes two cards and on each writes a different number.
Alice randomly chooses  one of these cards and peeks at the number.  
She must then decide which of the numbers is the larger: the number in her hand or the other, concealed number.
She wins each time her guess is correct.

By blind guessing Alice wins
with probability  $1/2$. Interestingly,
there is  a better strategy.
The idea is to compare the observed number $x$  with a 
threshold  number $t$,  sampled from some  probability distribution,  
for instance the normal.
If $x\geq t$ the strategy accepts  $x$ as the larger  number of two on the cards,
otherwise rejects $x$ and decides that the concealed number is the larger. 
By symmetry, Alice is right or wrong with the same probability, if $t$ falls either above
or below {\it both} Bob's numbers. But if $t$ falls {\it between} the two numbers,  
the decision is correct for sure:
the observed $x$ will be  accepted if it is the larger number and rejected 
otherwise. 
 The advantage of the strategy is  due to the latter possibility,  which has nonzero  probability 
whichever Bob's numbers happen to be.

For many  people this appears counter-intuitive, since
a random choice
from two  unknown numbers 
 seems   to bear   no  information on how the numbers  compare to one another.  
The history of the paradox can be traced back
to the work on theoretical statistics by
David Blackwell \cite{Blackwell} and Bruce Hill \cite{BHill}. 
Much of the interest was sparked  by Tom Cover's
\cite{Cover}   half-page abstract, where the game framework was explicitly introduced.
Steve Samuels \cite{Samuels} noticed 
that  the game belongs to the circle of questions around 
Martin Gardner's {\it Game of Googol} and 
 {\it Secretary Problem},
 where the guessing player aims to stop at the maximum of a sequence of unknown numbers observed  in random order  \cite{ BG,  Ferguson, SamSurvey}.  
The conundrum  was popularised  by 
Thomas Bruss \cite{Bruss} and
Ted Hill \cite{THill}, and 
many further fascinating  connections  with the
classical probability problems  were subsequently found \cite{2env, GGoogol, GermGoogol, GKrengel1, GKrengel2, HillKrengel, Samet, Winkler}.

From the perspective  of the zero-sum game theory the guessing game is fair, and
the blind guessing  is an equalising minimax strategy. 
The advantage of a threshold strategy over blind guessing
 is an interesting dominance phenomenon
that cannot occur in games with finitely many pure strategies. 
 Bob has no minimax strategy, and so the game has  no solution although has a definite value.
 The latter means that Bob is capable of  keeping Alice's advantage  arbitrarily small.

A common fallacy associated with the 
 game is that a smart strategy for Bob were choosing the numbers close to one another.
This  recipe is good to keep the advantage of a randomised threshold strategy small, but may fail
 when Alice uses some other strategy.
For instance, if  Bob employs a pure game strategy involving any two given numbers, then Alice can respond with
a sure-fire guessing.  Much more importantly, Bob should  really take care of  that the odds 
 given the observed number are about equal.
Such epsilon-minimax
strategies for the number-writing player have been designed in 
 \cite{ BG, Ferguson, GGoogol, GermGoogol}.

In this paper we 
discuss strategies for both players,  thus giving a complete account of the game.
In particular, we will characterise all Alice's minimax strategies, that guarantee the winning probability $1/2$.
We shall also  consider variations  of the game by granting the number-writing
player  different powers, 
to stress  the roles  of the admissible number range and
 exchangeability
 implied  by random choice of a card.
On this way
we  will  introduce 
a more involved game with two piles of cards.
In the new game Bob takes, let's say, 52 cards  
and on each writes a different number. The cards are turned face down, shuffled and dealt  in two 
piles of given sizes. 
Alice collects the first pile,  opens the cards and guesses if the number largest of all 52 numbers
is in her hand or in the second pile.   
This can be regarded as a version of the Game of Googol, where the stopping player has some return
options.
We will see that  properties of the two-pile game largely depend on whether the deck is divided evenly or not.
If   Alice gets exactly half of the deck,  she
can  take some advantage from learning the numbers by employing a threshold strategy, just like in the  game with two cards.
If the piles are of different size, Bob has a smart way to  generate the numbers, such  
that even if  Alice finds out his algorithm,
she  cannot do better than with the obvious strategy of betting   on the bigger pile. 

Although the paper is mostly self-contained, 
we shall assume the reader is familiar with  random variables and probability 
distributions on the level of a `second course' in probability or statistics.  
For fundamentals of the game-theoretic approach to statistical decision theory we refer to the classic text \cite{BlGi}.

\section{The general setting}

Before embarking on a discussion of the paradox
 it is helpful to consider a wider framework  for  the guessing games.
  Suppose that the revealed number $x$ and the concealed number $y$ are sample values of two random variables $X$ and $Y,$ respectively.
We assume   $P(X=Y)=0$, that probability of a `tie' is zero,  i.e. the numbers are different and so one of them is always  the larger and the other is the smaller. 
The {\it joint} probability distribution of the pair $X, Y$ is a strategy for Bob.
The rules of a particular game specify, sometimes implicitly,  a class of  joint distributions that Bob is  allowed to use.
The rules of the game belong to the common knowledge of players.

Having observed a value $x$, Alice is willing to guess which of the numbers is the larger. 
If she decides by tossing
a fair coin, she is right with probability $1/2$. Whether there are better strategies 
 depends  on the rules of the game.

  A pure strategy for Alice will be identified with a set $D\subset {\mathbb R}$ called a {\it decision set}.
  According to strategy $D$,
  the observed number $x$ is accepted as the larger if and only if  $x\in D$.
    To be able to define probabilities associated will the strategy we will assume that the set $D$ is Borel.

 Alice has two pure {\it blind guessing} strategies: one with 
  $D={\mathbb R}$ always accepting $x$, and the other with $D=\varnothing$ always rejecting $x$.
A pure {\it threshold strategy} has 
 $D=[t,\infty)$ with some fixed threshold $t$.

\vskip0.2cm
\noindent
{\bf Example A}  
It is easy to find a solution to the unconstrained game, where the joint distribution for $X,Y$ can be arbitrary.
Then playing $X=0$ and $Y=\pm 1$ with equal probabilities, Bob achieves that
Alice  wins with probability $1/2$ no matter how she decides.
Thus $1/2$ is the value of this game.

\vskip0.2cm
\noindent
{\bf Example B}
The game is more interesting when Bob is in some way restricted in choosing the joint distribution for $X,Y.$  
In the arrangement game introduced in \cite{GKrengel1} he 
 is {\it given} two numbers sampled independently from the uniform distribution on $[0,1]$,
but may choose which of the numbers is shown to Alice. 
The property characterising the relevant class of distributions is that the pair $X,Y$ arranged in increasing order
has uniform distribution on $\{(a,b): 0<a<b<1\}$. 
Still, Bob has
 an equalising strategy: show the number which is closer to $1/2$.
 By symmetry  of the uniform distribution about the median, the hidden number is equally likely to be bigger or smaller than the shown number. 

\vskip0.2cm

In line with the paradigm of  the zero-sum game theory, when choosing a strategy Bob should expect Alice's response most favourable for her.
The best response is an optimal strategy in a Bayesian decision problem, where
Alice knows exactly the joint distribution of $X, Y$ and can inambiguously evaluate strategies to find 
the one with the highest  winning probability.   
Having observed $x$, she will be able to
calculate the conditional probability
of winning  by accepting  the number in hand, that is 
$$\pi(x):=P(X>Y|X=x)=P(Y\leq x|X=x).$$
If she plays pure strategy $D$, the total probability of the correct guess is equal to
\begin{eqnarray*}
P(X>Y, X\in D)+ P(X<Y, X\in D^c)=\\E[\pi(X)1(X\in D)+(1-\pi(X))1(X\in D^c)],
\end{eqnarray*}
where $1(\cdots)$ equals 1 if $\cdots$ is true and equals 0 otherwise,
 $D^c$ is the complement of set $D$ and $E$  denotes the expectation.
The random variable  $\pi(X)$ is obtained by substituting $X$ in 
 the function $\pi:{\mathbb R}\to [0,1]$, so the
distribution of $\pi(X)$ is determined by the joint distribution of $X,Y.$
A minute reflection shows that the total probability of the correct guess  will be maximised if Alice accepts $x$ when $\pi(x)>1/2$,
and rejects $x$ when $\pi(x)<1/2$; while  decision in the case $\pi(x)=1/2$ can be arbitrary\footnote{In fact, optimality of this  strategy is a consequence 
of the principle of {\it dynamic programming}.}. 
Playing so, when
 $X=x$ she will be right with  probability
$$\max(\pi(x), 1- \pi(x))=\frac{1}{2}+\left|\pi(x)
-\frac{1}{2}\right|,$$
hence the  overall probability of the correct guess becomes
\begin{equation}\label{WinUnc}
v:=\frac{1}{2}+E\left|\pi(X)-\frac{1}{2}\right|.
\end{equation}
Note that the maximal probability $v$ is achieved by a pure strategy $D=\{x: \pi(x)\geq 1/2\}$. 
But if $\pi(X)$ assumes value $1/2$ with positive probability, there are also 
 mixed strategies, with randomised decision set, that achieve (\ref{WinUnc}).

A useful related concept is the  median of  the conditional distribution of $Y$ given that $X=x$ 
$$\mu(x)=\inf\left\{y: P(Y\leq y|X=x)\geq \frac{1}{2}\right\}.$$
In view of inequality
$$P(Y\leq\mu(x)|X=x)\geq \frac{1}{2}$$
the optimal decision set  can be determined by the condition 
\begin{equation}\label{viamu}
x\geq \mu(x).
\end{equation}

\vskip0.2cm
\noindent
{\bf Example C} Here is an important example. 
 Suppose $X$ and $Y$ are independent and identically distributed (iid) with some known continuous distribution. 
The Bayesian decision problem is then the simplest case of the
 {\it Full-Information Secretary Problem} \cite{SamSurvey}.
By independence
$\pi(x)=P(Y\leq x)=P(X\leq x)$ 
is the common cumulative distribution function, and the random variable $\pi(X)$ is uniformly distributed on $[0,1]$
(the fact known as the probability integral transform).
Since $|\pi(X)-1/2|$ is uniformly distributed on $[0,1/2]$, formula (\ref{WinUnc}) readily yields $v=3/4$. 
The conditional median is  a constant $\mu$ obtained as
the minimal  solution to the equation $\pi(t)=1/2$. Alice's
optimal strategy is of the threshold type, with threshold coinciding with the median.

\section{The paradox of two cards}

\subsection{Threshold strategies}
We turn now to  Tom Cover's game, where
Alice picks a card at random.
This has the effect of {\it exchangeability} \cite{Chow},
meaning that permutation $Y, X$ has the same joint distribution as $X,Y.$
In terms of the bivariate distribution function, exchangeability holds if 
$$P(X\leq x, Y\leq y)=P(Y\leq x, X\leq y)$$
for all $x,y$.
In particular, the random variables $X$ and $Y$ are identically distributed and $P(X>Y)=P(Y>X)=1/2$.
A pure strategy for Bob in this game can be modeled as a pair of fixed numbers $a,b$ with $a<b$, so that the associated exchangeable variables
$X,Y$ are these numbers arranged in random order.

We note in passing that
if $X$ and $Y$ are iid (as in Example C) then they are exchangeable,
but exchangeability  and independence are  very different properties.
Thus  for iid $X,Y$ the `no tie' condition $P(X=Y)=0$ means that $P(X=x)=0$ for every $x$, that is 
the distribution of $X$ is continuous.
In contrast to that,  exchangeability does not rule out discrete distributions, e.g. for each Bob's pure strategy the range of
$X$ has two elements.

 Suppose Alice employs a randomisation  device to generate a sample value $t$ of random variable $T$.
To exclude cooperation between the players, as inherent to the zero-sum game, 
 we shall assume that 
 $T$ is independent of $X,Y$.
 The associated threshold strategy accepts $x$  if $x\geq t$, and rejects   $x$ otherwise.
 By independence and exchangeability the triple $T,X,Y$ has the same joint distribution as  $T,Y, X$, thus 
  \begin{eqnarray*}
 P(X\geq T, Y\geq T, X>Y)=P(X\geq T, Y\geq T, Y>X) = \frac{1}{2}P(X\geq T, Y\geq T),\\
 P(X< T, Y< T, Y>X)= \frac{1}{2}P(X< T, Y < T),~~~~~~~~~~~~~~~~~~~~~~~~~~~~~~~\\
 P(X< T\leq Y)=P(Y< T\leq X)=\frac{1}{2}P(X< T\leq Y)+\frac{1}{2} P(Y<T\leq X).
 \end{eqnarray*} 
 Summing  probabilities on the right-hand side yields $1/2$, hence adding  to the sum $P(Y< T\leq X)$ the probability of correct guess with the threshold strategy becomes
 \begin{align}
w:=
\frac{1}{2} +P(Y< T\leq X). \label{WinP}
\end{align}

Strict inequality  $w>1/2$ 
 is ensured (for arbitrary distribution of $X,Y$) if $T$  is {\it fully supported}, i.e.  with  cumulative distribution function  
 \begin{equation}\label{cdf}
 F(t)=P(T\leq t)
 \end{equation}
strictly increasing for  $t\in {\mathbb R}$.
 The latter  condition is  equivalent to $$P(a<T\leq b))=F(b)-F(a)>0$$ for all $a<b$,  which is true whenever
$T$ has everywhere positive density.

Observe that
formula (\ref{WinP}) is valid under a weaker assumption that $T,X,Y$ have  the same joint distribution as $T,Y,X$
(that is, $X,Y$ are conditionally exchangeable given a value of $T$).
The case where all three variables $T,X,Y$ are exchangeable suggests the following  variation of the game,
in the spirit of Secretary Problem with a `training sample' \cite{SamTrS}.

\vskip0.2cm
\noindent
{\bf Example D} For readers unfamiliar with the Secretary Problem
this application  of (\ref{WinP})
might appear even more surprising than the paradox of two cards. 
Suppose Bob writes three distinct numbers on separate cards, which are shuffled and dealt face down in a row. 
Alice turns face up the first  and the second card, and must then guess 
 if the number on the 
second card is larger than a concealed number on the third card.
Using the number on the first card as a threshold to compare with the number on the second card, 
Alice wins with probability
$w=1/2+P(Y<T<X)=1/2+1/6=2/3$,
because due to exchangeability all $3!=6$ rankings of $T,X,Y$  are equally likely.
This example stresses the role of exchangeability in affecting the odds.
Although three numbers are a priori unknown, Alice can substantially benefit from comparing the first two.
\vskip0.2cm

For any given joint distribution of $X, Y$
the best of threshold strategies  has nonrandom threshold found by maximising $P(Y< t\leq X)$.
To make this probability positive, $t$ must lie 
between supremum and infinum of the range of $X$, although 
for some  such $t$ the probability can be zero\footnote{The assertion on top of p. 290 in \cite{Samuels} needs to be corrected.}.
For instance, if Bob plays a mixture
of pure strategies $1,2$ and $7,8$  we will have  $P(Y<t\leq X)=0$ for $t=5$.
If $X$ and $Y$ are not iid, typically  the optimal threshold does not coincide with  the median of  $X$.

\subsection{A probability inequality}

Once the first surprise is over, one might attempt explaining the paradox by arguing that 
the bigger the observed number $x$, the more likely it is the larger. 
However, exchangeability beats intuition also in this respect. 
Contrary to  the iid case, in general
the conditional probability $\pi(x)$ is neither  increasing in $x$,  nor even
exceeds $1/2$ for all sufficiently big $x$.
As a consequence, the optimal Alice's counter-strategy against $X,Y$ with given distribution may not be of the threshold type.

\vskip0.2cm\noindent
{\bf Example E}
Here is an example of erratic behaviour. For every  integer $j$ there is a pure strategy  for Bob with numbers $3j, 3j+1$.
Take a mixture of these pure strategies with some positive weights, that is pick $j$ from any probability distribution 
on the set of integers, then  arrange the numbers $3j$ and $3j+1$ in random order.  
For $X,Y$ thus defined, 
$\pi(x)=0$ for $x$ a multiple of $3$, but $\pi(x)= 1$ for $x=1$ modulo 3.
Obviously, Alice's best response is not a threshold strategy, rather 
a sure-fire strategy deciding that $x>y$ whenever $x=1$ modulo $3$, so with decision 
set $D=\{3j+1:\, j\in {\mathbb Z}\}$.
\vskip0.2cm

A similar example illustrates that $P(X>Y|X\geq t)$ need not be increasing in $t$. 
What {\it is} the correct intuition, is that 
the message `the observed number exceeds a threshold' increases the likelihood of  the event $X>Y$, because $X$  is certainly the larger when 
 the  hidden number falls below the threshold, whereas  the variables  are indistinguishable otherwise.
 Likewise, the message `the observed number falls below a threshold'
 speaks against  $X>Y$ 
because the hidden number may still exceed the threshold.
Formally, this conditional probability in focus is at least $1/2$, which is precisely the inequality underlying the paradox.
We will show this for a random threshold
$T$ independent of exchangeable $X$ and $Y$. Indeed,
\begin{align*}
P(X>Y,X\geq  T)&=P(X>Y,X\geq  T, Y\geq  T)+P(X>Y,X\geq T,Y< T) \\
&=
\frac{1}{2}P(X\geq T, Y\geq T)+P(Y< T\leq X)\\
&= \frac{1}{2} \left\{P(X\geq T, Y\geq  T)+ P(Y< T\leq X)\right\}   +\frac{1}{2}P(Y< T\leq X)\\
& =\frac{1}{2} P(X\geq  T)   +\frac{1}{2}P(Y< T\leq X).
\end{align*}
Assuming  that $P(X\geq T)>0$  we arrive at
\begin{equation}\label{wineq}
P(X>Y|X\geq T)=\frac{1}{2}+\frac{1}{2}P( Y< T|X\geq T)\geq \frac{1}{2}\,,
\end{equation}
where the bound is strict if $P(Y< T\leq X)>0$ . 

To conclude,  given the observed number exceeds a threshold,   betting  that this number is the larger is, overall, favourable.
However,  for a particular $x$, no matter how large, the bet may turn unfavourable.

Similarly, a negative outcome of the threshold comparison  reduces the likelihood
$$P(X>Y|X<T)=\frac{1}{2}-\frac{1}{2}P( Y\geq T|X < T)\leq\frac{1}{2}.$$
 This obvious counterpart of (\ref{wineq}) will be useful in what follows
to introduce deformations of threshold strategy.

There are several variations on the theme
how comparison with a threshold affects ranking.
One generalisation of  (\ref{wineq}) is the following.
 Let $X_1,\dots,X_n$ be exchangeable random variables without ties. The rank $R$ of $X_1$ is defined as the number of $X_j$'s not greater than $X_1$.
By exchangeability the distribution of $R$ is uniform on $n$ integers. 
Now, it can be shown that
for $T$ independent threshold, condition $X_1>T$ increases the probability of $R>k$ for every $k=1,\dots,n-1$.
That is to say,  the information that a number chosen at random out of $n$ numbers  exceeds a threshold stochastically increases the rank
of this number. See \cite{GKrengel2} for applications of relations akin to (\ref{wineq}) in optimal stopping.

\subsection{T. Hill's finite game}

In  the wonderful expository  article  \cite{THill} Ted Hill invites the reader
to solve a guessing game where  Bob is only allowed to use integer numbers from $1$ to, let us say, $m+1$.
This game reduces to a finite game, i.e. with finite sets of pure strategies for both players, because
for Alice it is sufficient to use decision sets $D\subset\{1,\dots,m+1\}$.
For the finite game the celebrated 
von Neumann's Minimax Theorem ensures that the game has a solution\footnote{Ted Hill provokes with a  fallacious argument:
`it  also seems obvious that the number-writer would never write a 1, since if [Alice turns] over a 1, [she] will always win by [not choosing 1]. But if he never writes a 1, he then would never write a 2 either since he never wrote a 1, and so on {\it ad absurdum}'.}.

Define a mixed strategy for Bob
by choosing $\beta$  uniformly from the set of integers
$1,2,\dots,m$. Given $\beta$,  the variables $X,Y$ are the numbers $\beta, \beta+1$ arranged in random order.
One can think of sampling a coin from a bag of $m$ coins with labels  $\beta, \beta+1$ on the sides, then tossing the coin  and noting a number $x$ on the upper side.
The resulting joint distribution of $X,Y$ is  uniform on the set of $2m$ integer pairs
\begin{equation}\label{m-distr}
\{(\beta, \beta+1):  1\leq \beta\leq m\}\cup \{(\beta+1, \beta):  1\leq \beta\leq m\}.
\end{equation}
When the observed number is  some $1< x< m+1$ both options $\beta=x$ and $\beta=x-1$ are equally likely, that is 
$$P(X=x)=P(X=x,Y=x+1)+ P(X=x,Y=x-1)=   \frac{1}{m}$$
and
$$\pi(x)= \frac{P(X=x,Y=x-1)}{P(X=x)}=    \frac{1}{2}.$$ 
For the extremes we have 
$$P(X=1)=P(X=m+1)=\frac{1}{2m}$$
and
$\pi(1)=0, \pi(m+1)=1$. Thus (\ref{WinUnc}) gives the total probability of win by the best response
$$v=\frac{1}{2}+ \frac{1}{2m}\left|0-\frac{1}{2}\right|+ (m-1)\cdot\frac{1}{m}\cdot0+  \frac{1}{2m}\left|1-\frac{1}{2}\right|=\frac{1}{2}+\frac{1}{2m}.$$

Define a threshold strategy for Alice with $T$ uniformly distributed over the set $\{2,\dots,m+1\}$.
Whichever distinct $x,y$ from $\{1,\dots,m+1\}$, the larger of them belongs to $\{2,\dots,m+1\}$.
Thus, whichever Bob's strategy
the event $T=\max(X,Y)$ is independent of $X,Y$ and has probability $1/m$. It follows that probability (\ref{WinP}) can be estimated as
\begin{eqnarray*}
w\geq \frac{1}{2}+P(Y<T=X)=   \frac{1}{2}+P(Y<X, T=\max(X,Y))&=& \\
   \frac{1}{2}+P(Y<X)P(T=\max(X,Y))&=&\frac{1}{2}+\frac{1}{2m}.
\end{eqnarray*}
We see that Alice can guarantee $1/2+1/(2m)$, while Bob has a strategy to bound her winning probability by this very value.
It follows that $1/2+1/(2m)$ is the value of the finite game, and that the described strategies of players are minimax.

Next, assuming that Bob can employ arbitrary integer numbers, 
we are in the situation where each player has infinitely many pure strategies,
hence the Minimax Theorem cannot guarantee existence of the value. 
However,
letting $m\to\infty$ in the finite game,  it is clear that Bob can keep Alice's winning probability as close to $1/2$ as desired.
 On the other hand,  since for the optimal Alice's counter-strategy $v\geq w>1/2$,
where $w$ is the winning probability for some threshold strategy,
Bob cannot achieve $1/2$ exactly. This implies that
the game with integer numbers has  value $1/2$.

But then  $1/2$ is also the value in the game where Bob
can use arbitrary real numbers. Bob has no minimax strategy and so the guessing game with real numbers
has no solution.

In contrast to the game on integers, restricting Bob's 
numbers to a finite interval will
not change the quantitative picture.
Indeed,  monotonic functions preserve the order relation hence  the game with arbitrary real numbers can be transformed 
in equivalent game where the range of $X,Y$ is bounded. The interval supporting 
$X,Y$ can be as small as wanted.

To disprove a common belief, we notice that Bob also has good strategies where $|X-Y|$ is big.
 For instance, he can play a mixture of
pure strategies $\beta k, \beta (k+1)$, with $\beta$ sampled from the uniform distribution on $\{1,\dots,m\}$ and $k$ fixed.
Such a strategy is epsilon-minimax for $k/m$ sufficiently small.
The latter does not exclude that $k=|X-Y|$ itself is a big number, nevertherless any
 random threshold $T$ will fall rarely between $X$ and $Y$.

\subsection{Minimax and admissible strategies}

Game-theoretically,
both pure 
blind-guessing  strategies with decision sets ${\mathbb R}$ and $\varnothing$, respectively,
are optimal (i.e. minimax) as they ensure
 the value of the game $1/2$. The same applies to any
 randomised blind guessing, e.g. deciding by tossing a fair coin.
However, every pure threshold strategy $D=[t,\infty)$ is minimax as well, and dominates  blind guessing.
A mixed threshold strategy with fully supported $T$  
strongly dominates the blind guessing, in the sense that $P({\rm win})>1/2$  no matter how Bob plays.
This  bizarre  dominance phenomenon cannot occur in finite games, where a strictly dominated strategy cannot be minimax.

In a discussion that came after the prize-winning  article \cite{Bruss},
Thomas Bruss  and Tom Cover raised the following natural questions: 
\begin{itemize}
\item[(i)]
Is the guessing strategy that always succeeds with probability 
strictly greater than $1/2$ a threshold strategy?
\item[(ii)] What is an optimal strategy? Is it of the threshold type?
\item[(iii)] What could be an alternative to threshold comparison, to benefit from learning one of the numbers?
\end{itemize}

Question (i) is answered in negative by means of a very simple construction.
 With every guessing strategy one can associate a {\it dual} strategy, which always makes the opposite decision. When the strategy wins the dual loses and vice versa. 
Taking a mixture of a threshold strategy and its dual, with weights, say, $3/4$ and $1/4$ (respectively), we obtain a strategy winning with probability  still
greater than $1/2$.
However, this  example is not very exciting, because the mixture is dominated by the threshold strategy.

To answer two other questions one needs to be  precise about the class of feasible  guessing strategies, and especially
about the very concept of randomised strategy.  
 Once the framework is established,
one should in the first instance focus on admissible (undominated) guessing strategies,  i.e. those 
which  cannot be improved in all situations, because the `optimal' strategies must be among admissible whatever `optimality' means.
A minimax guessing strategy will be called {\it superminimax} if  for any Bob's play the strategy  wins with probability strictly greater than $1/2$
and is admissible.   In response to the Bruss-Cover questions we will give  a complete characterisation of  the superminimax strategies.
In particular, we 
will   show that every admissible minimax strategy, in all what concerns its performance, is equivalent to a threshold strategy, and so from the optimisation viewpoint
   the class of threshold strategies is sufficient. With the reservation that there are many equivalent strategies,  
   the answer to (iii) is `no alternative'.

  \subsection{Characterising superminimax strategies}

Let us first recall some terminology.
For $\widetilde{D}$ and $D$  decision sets, we say that Alice's  pure strategy
$\widetilde{D}$ dominates $D$ if for every strategy of Bob
the probability of win with $\widetilde{D}$ is at least as great as with $D$,
and there exists Bob's strategy such that the probability is strictly greater. 
We say that $\widetilde{D}$ strongly dominates $D$ if the probability of win with $\widetilde{D}$ is always strictly greater than that with
$D$. A guessing strategy is called {\it admissible} if no other strategy dominates it.
For the purpose of checking  dominance it is sufficient to compare the outcomes
when Bob plays a pure strategy, because the general case follows by computing averages.
To verify admissibility of Alice's strategy one needs, in general, to make comparisons of outcomes also
under mixed strategies of Bob.

 The game in `normal form'  is specified by the outcomes when  the players use pure strategies.
If Alice plays $D\subset{\mathbb R}$ and Bob plays fixed $a<b$
the probability 
of correct guess is
\begin{equation}\label{ppay}
p(a,b;D):= \begin{cases}   \frac{1}{2}, ~~{\rm if}~a\in D,~b\in D,\\
\frac{1}{2}, ~~{\rm if}~a\in D^c,~b\in D^c,\\
1, ~~{\rm if}~a\in D^c,~b\in D,\\
0, ~~{\rm if}~a\in D,~b\in D^c.
\end{cases}
\end{equation}
where $1/2$ appears due to exchangeability of $X,Y$.

 Although simple, the analysis of pure strategies gives some hints.
 Identifying $D$ with its indicator function
 \begin{equation}\label{indic}
 F(x)=1(x\in D)
 \end{equation}
 the entry 0 in (\ref{ppay}) is the situation of {\it inversion} in $F$:
 \begin{eqnarray*}
 a<b, ~~F(a)>F(b).
 \end{eqnarray*}
 Inspecting (\ref{ppay}) it is clear that  
  $D$ is minimax if and only if   $F$  has  no inversions.
But such $D$ is either $D=\varnothing$ or ${\mathbb R}$ or  
 a right-tailed halfline, either $(t,\infty)$ or $[t,\infty)$ for some $t\in{\mathbb R}$.

 Regarding the admissibility of pure strategies, we note that
 if $D$ is bounded from the above, then this strategy is dominated by $\widetilde{D}=D\cup [\sup D,\infty)$, 
 because  $\widetilde{D}$ has less inversions, while other entries of (\ref{ppay}) are the same as for $D$. 
 Likewise, if $D^c$ is bounded from below, then $\widetilde{D}=D\setminus (-\infty,\inf D^c)$ dominates $D$.
 On the other hand,
 if $\sup D=\infty, \,\,\inf D^c=-\infty$ then $D$ cannot be dominated by a pure strategy.
Indeed,    suppose $\widetilde{D}$ dominates $D$.  For $a\in D^c, b\in D$  and $a<x<b$ it holds that if 
$x\in D$ then $1=p(a,x;D)\leq p(a,x,\widetilde{D})$ implies $p(a,x,\widetilde{D})=1$, hence $a\in \widetilde{D}^{\,c},\,x\in \widetilde{D}$.
If $x\in D^c$ then $1=p(x,b;D)\leq p(x,b,\widetilde{D})$ implies $p(x,b,\widetilde{D})=1$, and hence $x\in \widetilde{D}^{\,c}, b\in \widetilde{D}$. It follows that the sets  $D$ and $\widetilde{D}$ coincide between $a$ and $b$, and sending $a\to-\infty, \,b\to\infty$ along $D$ and  $D^c$, respectively,
we get $\widetilde{D}=D$.
It will be clear from what follows that this kind of argument still works when a mixed strategy is taken in the role of 
 $\widetilde{D}$, thus a pure strategy $D$ with 
$\sup D=\infty, \,\,\inf D^c=-\infty$ is admissible.

 We may now conclude that
 \begin{itemize}
 \item[] A pure guessing strategy is admissible and minimax if and only if the decision set is either 
  $[t,\infty)$ or $(t,\infty)$ for some  $t\in{\mathbb R}$.
 \end{itemize}

Turning to randomised strategies, 
statistical decision theory suggests two ways to introduce this concept:
\begin{itemize}
\item[(a)] by {\it mixing}, that is spreading a probability measure over the set of pure strategies,
\item[(b)]  by means of a {\it kernel}, which for every $x$ determines a distribution over the set of possible actions. 
\end{itemize}
The  relation between these two approaches has been thoroughly discussed, see e.g. \cite{Kirschner} and references therein.
In the context of the guessing game, the approach  (a)
 specialises as   a {\it random decision set} $\cal D\subset {\mathbb R}$.
To define a kernel  as in  (b) we need just one function
$F:{\mathbb R}\to[0,1]$  specifying probabilities $F(x)$ and $1-F(x)$ of 
two actions `accept $x$' and `reject $x$' (respectively) when the number $x$ is observed.
To link (a) to (b), with  random set $\cal D$ we 
associate the {\it coverage function}  \cite{Molchanov}
$$F(x):=P(x\in {\cal D}), ~x\in{\mathbb R},$$
which defines a kernel. We say that  two random sets are  {\it equivalent} if they have the same coverage function. 
Note that the idea of equivalence class leaves open alternative concepts of a strategy associated with given $F$.

If Bob can only use integer numbers, then with every $F:{\mathbb Z}\to[0,1]$ we can associate a random set by means of  a Bernoulli process, that is  by including $x$ in $\cal D$ with probability $F(x)$, independently for all 
$x\in{\mathbb Z}$. For the game on reals the analogous construction will lead to a set which is not measurable, hence the connection between (a) and (b) is more subtle.

To avoid  annoying complications we shall assume that $F$ is {\it cadlag}, that is  right-continuous with left limits.
The instance of randomised threshold strategy neatly fits in the random set framework: in this 
case  ${\cal D}=[T,\infty)$, and  since $x\in {\cal D}$ means the same 
 as $T\leq x$  the coverage function is the cumulative  distribution function (\ref{cdf}) of  the random threshold $T$.
 To avoid separate  treatment of   blind-guessing strategies it will be convenient to 
consider them as  special cases  of the threshold strategy,  with $-\infty$ and $\infty$ being legitimate values for the threshold. A threshold strategy with $P(T=\pm\infty)=0$ will be called {\it proper}, hence satisfying
 \begin{equation}\label{proper}
 \lim_{x\to-\infty} F(x)=0, ~\lim_{x\to\infty} F(x)=1.
 \end{equation}

When Bob plays exchangeable $X,Y$ with some given joint distribution and Alice uses $\cal D$ with some coverage function $F$ the outcome of the game is assessed as
\begin{eqnarray}
P({\rm win)}=E\,[1(X>Y)F(X)+1(Y>X)(1-F(X))]= \nonumber \\
E[1(X>Y)F(X)+1(Y>X)-1(Y>X)F(X)]=    \nonumber \\
\frac{1}{2}+ E\,[1(X>Y)F(X)-1(X>Y)F(Y)]=  \nonumber \\
\frac{1}{2}+ E\,[1(X>Y)\{F(X\vee Y)-F(X\wedge Y)\}]=  \nonumber \\
\frac{1}{2}+\frac{1}{2}\,E\{F(X\vee Y)-F(X\wedge Y)\}\label{prowin},
\end{eqnarray}
where $\vee, \wedge$ denote maximum and minimum, respectively.
The passage to the third line is justified by exchangeability.
The last line follows by the following observation: for symmetric function $\varphi$ from $1(X>Y)+1(Y>X)=1$ and exchangeability we have 
$$E[(1(X>Y)\varphi(X,Y)]=E[(1(Y>X)\varphi(X,Y)]=\frac{1}{2} E[\varphi(X,Y)].$$

 At a first glance, 
  formula  (\ref{prowin}) may seem to disagree with (\ref{ppay}), as the averaging will apparently involve more complex two-point coverage probabilities $P(a\in {\cal D}, b\in{\cal D})$. However,
re-writing (\ref{ppay}) as
\begin{eqnarray*}
p(a,b;D)= \frac{1}{2}+\frac{1}{2}\,1(a\in D^c, b\in D)-\frac{1}{2}\cdot\,1(a\in D, b\in D^c)=\\
\frac{1}{2}+\frac{1}{2}\cdot\,[1(b\in D)-1(a\in D, b\in D)]-\frac{1}{2}\,[1(a\in D)-1(a\in D, b\in D)]=\\
\frac{1}{2} +\frac{1}{2}\,[1(b\in D)-1(a\in D)],
\end{eqnarray*}
and averaging over the realisations of the random decision set, for the winning probability with  $\cal D$ against pure strategy $a<b$ we obtain
\begin{equation}\label{PPP}
P({\rm win})=\frac{1}{2}+\frac{1}{2}(F(b)-F(a))
\end{equation} 
in accord with (\ref{prowin}).  In particular,  for threshold strategy (\ref{prowin})   becomes (\ref{WinP}). 
 By the analogy with the threshold case, we call general $\cal D$  {\it proper}, if the coverage function satisfies 
 (\ref{proper}).

 It is very fortunate that probability (\ref{PPP}) is expressible via such a simple characteristics of $\cal D$. 
 It is immediate from the formula that
 \begin{itemize}
 \item[(i-a)] A mixed strategy $\cal D$ is minimax if and only if  the coverage function $F$ is nondecreasing,
 \item[(i-b)] A mixed strategy $\cal D$ strongly dominates blind-guessing 
 if and only if the coverage function $F$ is strictly increasing.
 \end{itemize}
 
 \noindent
 In particular,   every  threshold strategy with fully supported $T$, possibly with atoms at $\pm\infty$
 (blind-guessing component), is winning with probability strictly greater than $1/2$.  A more interesting example, stepping away from threshold strategies
 is the following.
 \vskip0.2cm
\noindent 
 {\bf Example F}
Let $T$ be a  fully supported random variable with (strictly increasing) cumulative distribution function $F$, possibly not proper.
 With $T$ we associate both the threshold strategy, and  a { dual} strategy which accepts $x$ when $x$ is {\it smaller} than the sample value of $T$.  The decision set for the dual strategy is $(-\infty,T)$, and 
 the coverage function is $1-F$.
 For $0\leq \gamma<1$
 define a {\it $(F,\gamma)$-strategy}  as a mixture of the threshold strategy and its dual with weights $\gamma$ and  $1-\gamma$.
 The coverage function for the $(F,\gamma)$-strategy is the convex combination
 $$F_\gamma(x)=\gamma F(x)+(1-\gamma)(1-F(x)).$$
 Clearly, the  $(F,\gamma)$-strategy is not of
 threshold type because with probability $1-\gamma$ the decision set is a left-tailed half-line.
 We have then for $a<b$
 $$P({\rm win})=
  \frac{1}{2}+\frac{1}{2}(F_\gamma(b)-F_\gamma(a))=
  \frac{1}{2}+\frac{2\gamma-1}{2}(F(b)-F(a)),$$
 which exceeds $1/2$ for $\gamma>1/2$.

  \vskip0.2cm
  
 We turn  to  the admissibility of mixed guessing strategies. A simple sufficient condition is 
 \begin{itemize}
 \item[(ii-a)] If $\cal D$ has a proper coverage function $F$ then $\cal D$ is admissible.
 \end{itemize}
 Indeed, suppose $F$ is proper and $\widetilde{F}$ is a coverage function for some strategy, such that
 $\widetilde{F}(b)-\widetilde{F}(a)\geq F(b)-F(a)$ for all $a<b$. Sending $a\to-\infty, b\to\infty$ we get 
 $\widetilde{F}(b)-\widetilde{F}(a)\to 1$, which is only possible when $\widetilde{F}$ is also proper.
 Therefore
 letting $a\to-\infty$ yields $\widetilde{F}(b)\geq F(b)$,
 and letting $b\to\infty$ yields $1- \widetilde{F}(a)\geq 1-{F}(a)$. Since $a<b$ arbitrary we have $F=\widetilde{F}$, hence $\widetilde{F}$ 
 cannot dominate $F$.

 Without engaging in fuller analysis of admissibility, as a partial converse to (ii-a)
 we have
 
 \begin{itemize}
 \item[(ii-b)] If $\cal D$ has  coverage function $F$  with $\liminf\limits_{x\to-\infty} F(x)>0$ 
 or $\limsup\limits_{x\to\infty} F(x)<1$  (so, not proper) then $\cal D$ is dominated. 
 \end{itemize}
 Indeed, the increments will be only bigger, if  we increase $F$  on some halfline $[x,\infty)$  by a constant,
 or decrease  on some $(-\infty,x)$ by a constant.

 If $F$ is not proper and nondecreasing (but not constant), we can dominate $\cal D$ by
 any $\widetilde{\cal D}$ with coverage function
 $$\widetilde{F}(x)=\frac{F(x)-c_1}{c_2-c_1},$$
 where $c_1=\lim_{x\to-\infty} F(x)$ and  $c_2=\lim_{x\to\infty}F(x)$.
For instance, the $(F,\gamma)$-strategy from Example F is dominated by the threshold strategy with 
 $T$ distributed according to $F$;  the dominance is strong if $T$ is fully supported.

If the coverage function $F$ is constant, then $\cal D$ performs like a blind-guessing strategy,
hence $\cal D$ is not admissible. Constant $F$ appears as a coverage function for a {\it stationary} (translation invariant) random set $\cal D$.

\vskip0.2cm
\noindent
{\bf Example G} Let $\cal P$ be the random set of atoms of a homogeneous Poisson point process with rate $\lambda$.
A realisation of the process is a random scatter of  isolated points (atoms). 
Define a decision set as the pointwise sum of sets ${\cal D}={\cal P}+[0,1]$, which is the union of unit intervals having left endpoints in $\cal P$. 
A point $x$ is uncovered by $\cal D$ with
 probability $e^{-\lambda}$, as this occurs whenever
$\cal P$ has no atoms in $[x-1,x]$, 
The coverage function is therefore constant $F(x)=1-e^{-\lambda}$,
hence this strategy always wins with probability $1/2$.

Instead of building upon a Poisson process, one can use any stationary renewal point process, expressing $F$  via the  distribution of  the age  (also known as backward excess) random variable.

\vskip0.2cm

 Every  cadlag nondecreasing $F$ with values in $[0,1]$ is a distribution function for some random variable, hence  the  coverage function for a threshold strategy.  
Putting this together  with (i-b), (ii-a) and (ii-b) we arrive at a crucial conclusion

 \begin{itemize}
 \item[(ii-c)] A guessing strategy $\cal D$ is superminimax if and only if the coverage function $F$ is strictly increasing and proper.
 In this case $\cal D$ is equivalent to the threshold strategy with $T$ having the cumulative distribution function $F$.
 \end{itemize}

 We give a construction of  superminimax strategies not of the threshold type.

 \vskip0.2cm
 \noindent {\bf Example H} For $\lambda$ a positive locally integrable function on $\mathbb R$,  consider a (inhomogeneous) Poisson point process $\cal P$ of intensity $\lambda$.
Define ${\cal D}={\cal P}+[0,1]$.
 With positive probability the set $\cal D$ has many connected components, hence this is not a threshold strategy.
The coverage function is computed as 
 $F(x)=1-e^{-\Lambda(x)},$
where
 $\Lambda(x)=\int_{x-1}^x \lambda(t) dt$
 is the mean number of Poisson atoms in $[x-1,x]$. To get a proper strictly increasing $F$, it is enough to require that $\lambda$ be increasing, with limits at $\pm\infty$ being $0$ and $\infty$, respectively.
 To meet these conditions we can take the intensity function  $\lambda(t)=e^t$.
 \vskip0.2cm

Finally, we comment on the game on $\mathbb Z$.
If Bob is restricted to use only integer numbers the above arguments simplify in that there are no longer measurability and continuity concerns. A random set ${\cal D}\subset {\mathbb Z}$ can be identified with a $0$-$1$ process on the lattice. The simplest such process is homogeneous Bernoulli, which is equivalent to the blind guessing.
Borrowing an idea from 
the interacting particle systems  \cite{Liggett} and the theory of generalised exchangeability \cite{GO} a superminimax strategy can be
constructed as a $q$-deformation of the Bernoulli process, as in the following example. 
\vskip0.2cm

\noindent
{\bf Example I} There exists a random bijection $\Sigma:{\mathbb Z}\to{\mathbb Z}$ with the following property.
For every $n$, and a subset $S\subset{\mathbb Z}$ with $2n+1$ elements, 
given that the (unordered) set of values of 
$\Sigma|_{[-n,n]}$  is $S$, with probability proportional to $q^I$
the bijection $\Sigma|_{[-n,n]}\to S$
coincides with any of  $(2n+1)!$ bijections  $\sigma:{[-n,n]}\to S$ where $I$ is the number of inversions in $\sigma$. 
We can view $\Sigma$ as an infinite analogue of the Mallows model for finite random permutations, see \cite{GO} for details.
The decision set will be defined as ${\cal D}=\{j: \Sigma(j)\geq 0\}$. The indicator function of $\cal D$ may be thought of as configuration  of `particles' occupying some 
positions on the lattice. This configuration is obtained by using $\Sigma$ to shuffle the initial configuration where particles occupy positions $j\geq 0$ while
 positions $j<0$ are the `holes'.
For $0<q<1$ the bijection $\Sigma$ tends to be increasing, and   $\cal D$  has the decisive property 
that the pattern `hole, particle' in positions $a<b$ is more likely than `particle, hole'; thus $\cal D$ is a superminimax strategy.
The random configuration $\cal D$ is an equilibrium state  (also known as `blocking measure') for asymmetric simple exclusion process \cite{Liggett}.


\subsection{Repeated game}

The objective interpretation of probability  involves the long-run frequency in a series of identical trials. 
Thus speaking of probability of an outcome
in the guessing 
 game we actually  that the game is played in many rounds, and that
 each player employs the same randomisation device in every round.
 In the time perspective, the latter means that the strategies are {\it stationary}.
Considering the guessing game as a {\it repeated game} \cite{Sorin}, we may allow playing
nonstationary strategies, by changing randomisation method from round to round.
The long run frequency of wins is a natural objective for the guessing player.
A solution to repeated game exists in stationary strategies  if the basic (stage) zero-sum game has a solution.
But in  the guessing game Bob has only epsilon-minimax strategies, therefore
 it perfectly makes sense to use nonstationary strategies for writing the numbers, in each successive round making guessing harder.
For instance,  Bob  will  bound the long-run frequency of Alice's wins by {\it exactly} $1/2$ by employing the uniform distribution on the finite set (\ref{m-distr}) in round $m$.
Then, in terms of the long run frequency,  the strategy of blind guessing cannot be  outperformed, and the repeated game has a solution.

\subsection{B. Hill's assumption and the invariance principle}

An ideal strategy for Bob is exchangeable pair $X, Y$ without ties satisfying
 $P(\pi(X)=1/2)=1$. This condition is the $n=1$ instance of 
 {\it assumption  $A_n$}  introduced  by  Bruce Hill \cite{BHill}  
 as a model for nonparametric predictive inference.
He used threshold comparison to show that  $A_1$ and then, by induction, $A_n$  is  impossible.
But Lane and Sudderth \cite{Lane} showed that  $A_n$ can be realised in the setting of finitely additive probability theory,
which leads to a kind of unconventional  minimax strategy for Bob.

The uniform distribution on (\ref{m-distr}) with large $m$ can be regarded as approximation to $A_1$. 
It is instructive to represent this strategy as
\begin{equation}\label{disc}
X={B}+G_1, ~Y=B+G_2,
\end{equation}
where the pair  $G_1,G_2$ is exchangeable, and independent of the
  {\it location parameter} $B$.
Using the additive form many other approximations to $A_1$, hence good strategies for Bob can be constructed. 
More insight comes from the
statistical Principle of Invariance which asserts existence of {\it invariant} minimax strategy in decision problem with certain group of symmetries.

\vskip0.2cm\noindent
{\bf Example J} Consider the game on $\mathbb Z$  with a  further constraint that Bob must use two consequtive integers.
A pure strategy for Bob amounts to choosing the location parameter $\beta\in{\mathbb Z}$, to play  exchangeable $X=\beta+G_1, Y=\beta+G_2$,
where $G_1,G_2$ is the pair $0,1$ arranged in random order.
The group $\mathbb Z$ acts in this setting by translations and is {\it admissible} (cf \cite{Blackwell}, Definition 8.6.2).
The latter means that, for all $\beta, x\in{\mathbb Z}$
when Bob plays pure strategy $\beta$ and Alice plays a strategy which accepts a given $x$ as the larger number when $x$ is observed,
the outcome is the same as in the game where
Bob plays $\beta+z$ and Alice  accepts $x+z$ when 
$x+z$ is observed. Since $\mathbb Z$ acts transitively on the range of observed variable,  any invariant minimax strategy must be constant,
hence it is blind guessing. 
Uniform distribution on $\{-m,\dots,m\}$ (for big $m$) is Bob's  epsilon-minimax strategy,  which may be regarded  as  approximation
to the ideal strategy for Bob, the latter being the infinite invariant measure on $\mathbb Z$.

\vskip0.2cm
\noindent
{\bf Example E (continued)}
Suppose $G_1,G_2$ in (\ref{disc}) is  the pair of numbers $0,1$ arranged in random order, and that $B$ takes values in the group $3{\mathbb Z}$ of integer  multiples of $3$.  
Alice's invariant minimax strategy is constant on the orbits $\{3j:j\in {\mathbb Z}\}$ and $\{3j+1:j\in {\mathbb Z}\}$, and it is the sure-fire strategy which 
accepts $x=3j+1$  and rejects $x=3j$.

\vskip0.2cm

\vskip0.2cm
\noindent
{\bf Example K} 
To approximate $A_1$ by continuous distributions  we can take 
\begin{equation}\label{loc}
X=B+U_1, Y=B+U_2\,,
\end{equation}
where $B$ is   uniformly distributed on $[-m,m]$  and independent of   $U_1,U_2$, and the variables $U_1,U_2$  are iid uniform on $[0,1]$.
Similarly to (\ref{disc}),  taking $m$ big  we  achieve that $\pi(X)=1/2$ with high probability, and
this yields a epsilon-minimax strategy for Bob in the game where he controls a continuous  location parameter $\beta\in {\mathbb R}$. 
The value of this game is $1/2$.
The advantage of using (\ref{loc}) , as  compared to discrete distributions, is that 
$X,Y$ is a mixture of iid sequences, hence can be extended to an infinite sequence $B+U_i$ of exchangeable random variables.
A version of the Game of Googol for such location mixtures was studied in \cite{Petruccelli}.

\vskip0.2cm
\noindent
{\bf Example L} 
Yet another approximation to $A_1$ is obtained 
via the multiplicative form
\begin{equation}\label{scale}
X=A U_1, Y=A U_2,
\end{equation}
where $A$ is a positive {\it scaling parameter} independent of uniform iid $U_1,U_2$. 
This can be interpreted in context of the game where Bob only controls the scaling parameter.
To have $E|\pi(X)-1/2|$ small, one can take $A$ with density 
$$P(A\in [a, a+da])= \frac{1}{2\log m} \frac{da}{a}\,,~~~m^{-1}<a<m,$$
which may be seen as a finite approximation to the invariant measure on the multiplicative group ${\mathbb R}_+$.

Scale mixtures of uniform distributions, based on another approximation to the invariant measure on  ${\mathbb R}_+$, will play important role in the remaining part of the paper.

\section{The two piles game}

\subsection{Setup}

We focus next on  a game that bears features of  asymmetric two-cards game, the  Game of Googol with $n$ cards
and the game with two sequences from \cite{Sz}.
Bob takes $n$ separate cards and on each writes a different number.
The deck is shuffled and dealt face down in  two piles of $k$ and $n-k$ cards. The pile sizes are known to both players.
Alice collects the pile with $k$ cards, looks at the 
numbers and guesses if the number largest of all $n$ numbers is in her hand or in the other pile with $n-k$ cards. 
She wins  if her guess is right.

We model Bob's strategy  as  a joint probability distribution for $n$  
random variables $Z_1,\dots,Z_n$ which do not tie, i.e. $P(Z_i =Z_j)=0$  for $i\neq j$, and are  exchangeable. 
Exchangeability is equivalent to the property that the sequences $\ldots, Z_i, Z_{i+1},\ldots$ and $\ldots, Z_{i+1}, Z_{i},\ldots$ 
with just two neighbours swapped have the same  joint distribution.

The numbers $z_1,\dots,z_k$ in Alice's hand will be identified with values of random variables   $Z_1,\dots,Z_k$, and the hidden numbers in the second pile with  values of
$Z_{k+1},\dots,Z_n$. Using $\vee$ for `maximum' let  $X=Z_1\vee\ldots\vee Z_k$ and $Y=Z_{k+1}\vee\ldots\vee Z_n$ be the respective maxima in piles. 
Having observed $z_1,\dots,z_k$, Alice's dilemma amounts to guessing which of the inequalities $x>y$ and $x<y$  is true for the values $x,y$ of $X,Y$, respectively.

The setting is similar to the game with two cards, but there are  two delicate distinctions.
Firstly, Bob determines the joint distribution of $X,Y$ indirectly,  via the joint distribution of $Z_1,\dots,Z_n$.
Secondly,  Alice has more information when observing $Z_1,\dots,Z_k$ as compared to just learning the maximum of these $k$ variables. 
Her best-response strategy is to decide that $x>y$ if  
\begin{equation}\label{crit}
P(X>Y|Z_1=z_1,\dots, Z_k=z_k)\geq  \frac{1}{2}\,,
\end{equation}
where, in general,  the  conditional probability depends on all $z_1,\dots,z_k$ in a complex way.
For tractability reasons it will be useful  to find a framework where this probability is expressible in terms of a few summary statistics of the observed numbers.

By exchangeability, 
$$P(X>Y)=\frac{k}{n}\,,$$
because each $Z_j$ is equally likely to be the largest of all $n$ numbers.
We call the game {\it symmetric} if the deck is split evenly  ($n$ even, $k=n/2$).

If the game is symmetric, $X$ and $Y$ are themselve exchangeable. Like in the two-cards game,  a threshold strategy dominates the blind guessing.
One might anticipate that the value of the game is still $1/2$, but the tools we used to show this in the two-cards game are not appropriate for $n>2$.

In the asymmetric game the  pile of bigger size is more likely to contain the overall largest  number. Under  `blind guessing'  we shall understand Alice's strategy which always decides on the larger pile, regardless of the revealed numbers.
Since  on the average $E[P(X>Y|Z_1,\ldots,Z_k)]=k/n\neq 1/2$,   we may hope   
to find a joint distribution for $Z_1,\dots,Z_n$, 
such that  the conditional probability
is always within $\epsilon$ from the constant $k/n$, hence   
 compares to $1/2$  in the same way as $k/n$ does. 
 This would also imply the same kind of definite relation between $X$ and the conditional median of $Y$ given $Z_1,\dots,Z_k$.
 For  such Bob's strategy  the blind guessing would be Alice's best response.

\subsection{The game with independent numbers}

Suppose Bob plays iid $Z_1,\dots, Z_n$   with uniform distribution on $[0,1]$, and Alice plays her best-response (Bayesian) strategy.
For $x=z_1\vee\ldots\vee z_k$, condition (\ref{crit}) becomes $P(Y< x)\geq 1/2$, which holds exactly for   $x\geq \mu$, where threshold $\mu$ is the median of $Y$.

From $P(X< x)=P(Z_1<x,\ldots,Z_k<x)=x^k$ the density function of 
$X$ is $kx^{k-1}$ for $x\in[0,1]$.  
From $$\frac{1}{2}=P(Y< \mu)=P(Z_{k+1}\leq \mu,\dots,Z_{n}\leq \mu)=\mu^{n-k}$$  the median  is $\mu=2^{-1/(n-k)}$.
The optimal probability of win is computed as
\begin{eqnarray*}
v&=&P(X>Y,X>\mu)+P(X<Y, X<\mu)\\
&= &\int_\mu^1 x^{n-k} \cdot kx^{k-1}dx+ \int_0^\mu (1-x^{n-k})kx^{k-1}dx \\
&=&\frac{k}{n} +\mu^k+ \frac{2k}{n}\mu^n\,,
\end{eqnarray*}
which upon substituting the median gives
\begin{equation}\label{WinInd} 
v=\frac{k}{n}+\left( 1-\frac{k}{n}   \right) 2^{-k/(n-k)}\,.
\end{equation}

For symmetric games $k/n=1/2$, and we arrive at $v=3/4$ in accord with our previous computation for the two-cards game.

Observe that $v$ is a function of  the ratio $k/n$. This is easy to explain using the invariance of the game under monotonic transformations of the number scale.
 For increasing function $\varphi: (0,1)\to{\mathbb R}$, the change of variables $\varphi(Z_1),\dots,\varphi(Z_n)$
yields an equivalent game, since $\varphi$ preserves the order relations. Arbitrary continuous distribution for  the random variable $\varphi(Z_i)$  can be obtained from uniformly distributed  $Z_i$.  
In the iid case, the  game with  pile sizes 
$\ell k, \ell(n-k)$ is  reducible to the game with pile sizes $k, n-k$ by passing to $n$ iid variables $Z_i'=Z_{\ell (i-1)+1}\vee\ldots \vee Z_{\ell i}$ ($i=1,\dots,n$)
and noting that $Z_1\vee\ldots \vee Z_{\ell k}=   Z_1'\vee\ldots\vee Z_k'$ and similarly for the other pile.

Minimising $v$ as function of $k/n$, the ratio worst for Alice is about $0.587$, with minimum $v$ being close to $0.741$.
On the other hand, the maximum $1$ is achieved for $k/n\to 0$  or $1$.

\subsection{A class of multivariate densities}

Multivariate distributions well suited for our purposes have densities the form 
\begin{equation}\label{diamond}
f_n(z_1,\dots,z_n)=g_n(z_1\vee\ldots\vee z_n), 
\end{equation}
where $z_j>0$ for $j=1,\dots,n$ and $g_n$ is a nonnegative function. 
Random variables $Z_1,\dots,Z_n$ with joint density $f_n$ are exchangeable because the maximum $z_1\vee\ldots\vee z_n$ is a symmetric function.
The total probability integral must be one,  therefore
$g_n$ must  satisfy 
\begin{align}
\int_0^\infty \dots \int_0^\infty f_n(z_1,\dots,z_n)dz_1\dots dz_n&=&   \nonumber  \\
n\int_0^\infty g_n(z_1)dz_1 \int_0^{z_1} \dots \int_0^{z_1} dz_2\dots dz_n&=& \nonumber \\
n \int_0^\infty g_n(t)t^{n-1}dt&=1.\label{inte}
\end{align}

 For $k<n$ the joint density  of  
$Z_1,\dots,Z_k$ 
is obtained by integrating out  the variables $z_{k+1},\dots,z_n$ as

\begin{eqnarray}\nonumber
f_k(z_1,\dots,z_k)=\int_0^\infty\dots\int_0^\infty g_n(x\vee z_{k+1}\vee \cdots\vee z_n) dz_{k+1}\dots dz_n=\\
\label{gk}
x^{n-k} g_n(x)+(n-k)\int_x^\infty y^{n-k-1}g_n(y)dy= :g_k(x),
\end{eqnarray}
where $g_k$, similarly to  $g_n$, is  a function of   maximum $x=z_1\vee\ldots\vee z_k$.
The conditional joint density of  $Z_{k+1},\ldots,Z_n$  at  locations $z_{k+1},\ldots,z_n$  
given the observed numbers $Z_1=z_1,\ldots,Z_k=z_k$ becomes $g_n(x\vee y)/g_k(x)$, where
$y=z_{k+1}\vee\ldots\vee z_n$.  
Integrating  out  
 $z_{k+1},\ldots,z_n$, with each variable running from  $0$ to $x$, we arrive at 
\begin{eqnarray}
\nonumber
P(X>Y|Z_1=z_1,\dots,Z_k=z_k)&=&
P(Y<x |X=x)=\\ 
\label{heartsuit}
\frac{x^{n-k} g_n(x)}{g_k(x)}&=:&\pi_{n,k}(x).
\end{eqnarray}

The rationale behind (\ref{diamond}) can be now appreciated.
The conditional density of $Y$ given $Z_1=z_1,\dots,Z_k=z_k$ depends only on the maximum $x$ of the observed numbers.
In this sense we are back to the two-cards problem with certain class of  bivariate distributions for $X,Y$.

\subsection{Scale mixtures of uniform distributions}

We want to find a function  $g_n$ for which $\pi_{n,k}(x)$ is close to the constant $k/n$ for all $x>0$.
An ideal candidate is  $g_n(x)=n x^{-n}$ as (\ref{heartsuit}) gives  $\pi_{n,k}(x)=k/n$ exactly. But this function 
cannot be used to define a proper density  (\ref{diamond}) since the integral in
 (\ref{inte}) diverges. In what follows we will find a proper  deformation of  the power function.

An important  subclass of multivariate densities (\ref{diamond}) is the family of scale mixtures of uniform distributions, which we now introduce.
Let $U_1, U_2,\dots$ be iid uniform on $[0,1]$, and let $A$ be a positive random variable, independent of $U_j$'s. 
Consider an infinite exchangeable sequence 
\begin{equation}\label{scale}
Z_j=AU_j, ~~j=1,2,\dots
\end{equation}
Suppose $A$ has  density $h$ on the positive halfline. Given $A=a$, the variables  $Z_1,\dots,Z_n$ are iid, uniform on $[0,a]$, hence their (unconditional) joint density is of the form (\ref{diamond}) with 
\begin{equation}
\label{gtoh}  
g_n(x)=\int_x^\infty a^{-n} h(a) da, ~~~n=1,2,\dots
\end{equation}
It is not hard to check that this formula agrees with (\ref{gk}) for $k<n$.
Conversely, if  $g_n$ is a nonnegative function satisfying (\ref{inte}) and  $g_n'\leq 0$ then $f_n$ is the $n$-dimensional density for (\ref{scale}) with variable $A$ having density $h(a)=-g'_n(a)a^n,~ a>0$.
See more in  \cite{GnedinSPL} on this connection.

Observe that for the improper distribution with density $h(a)=a^{-1}, a>0$ formula (\ref{gtoh}) yields the ideal function 
 $g_n(x)=nx^{-n}$.

To achieve convergence we take a piecewise smooth density
 \begin{eqnarray}\label{phi}
h(a)=
\begin{cases} 
c\,a^{-\delta-1}, ~~{\rm for~~}a\geq 1,\\
c\, a^{\delta-1}, ~~~~{\rm for~~}0<a<1,
\end{cases}
\end{eqnarray} 
with $0<\delta<1$ yet to be chosen and $c=2\delta^{-1}$. Then the integral (\ref{gtoh}) is easily calculated as

\begin{eqnarray}\label{gnon}
g_n(x)= 
\begin{cases}
\frac{c}{n+\delta}\, x^{-n-\delta}\,, ~~~~~~~~~~~~~{\rm for~~}x\geq 1,\\

 \frac{c}{n-\delta} \, x^{-n+\delta} -\frac{2c\delta}{n^2-\delta^2} , ~~~{\rm for~~}0<x<1,\\
\end{cases}
\end{eqnarray}
where for the second part we used the decomposition
$$g_n(x)=\int_x^1 a^{-n} (c a^{\delta-1})da +g_n(1).$$
Now   a part of (\ref{heartsuit}) is straight from the first part of (\ref{gnon})
\begin{equation}\label{p1}
\pi_{n,k}(x)=
 \frac{k+\delta}{n+\delta}, ~~~{\rm for}~x\geq 1,
\end{equation}
and similarly with slightly more effort  the other part becomes 
\begin{eqnarray}\label{p2}
\pi_{n,k}(x)=
\frac{k- \delta\left(1+  \frac{2(k-\delta)}{n+\delta}\, x^{k-\delta}\right)}    { n-\delta\left(1 + \frac{2(n-\delta)}{k+\delta}\, x^{k-\delta}\right)}, ~~~{\rm for~}0<x<1.
\end{eqnarray}
From (\ref{p1}) and (\ref{p2}) we see that for every  integer $N$ and given positive $\epsilon<1/N$ it is possible to choose $\delta>0$ so small that 
$$\left|\pi_{n,k}(x)-\frac{k}{n}\right|< \epsilon$$
holds
uniformly in $x> 0$ and $1\leq k\leq n\leq  N$.  With this choice  of parameter   
$\pi_{n,k}$ relates to $1/2$ exactly as $k/n$ does, unless $k=n/2$.
But then the distribution of $Z_1,\dots,Z_n$  defined by (\ref{diamond}) and  (\ref{gnon}) has the properties that

\begin{itemize}
\item[(i)]  for $k>n/2$ the optimal Alice's counter-strategy is to always accept the maximum of the first pile $x$,
\item[(ii)]  for $k<n/2$ the optimal Alice's counter-strategy is to always reject $x$,
\item[(iii)] for  $k=n/2$ the optimal Alice's counter-strategy wins with probability at most $1/2+\epsilon$.
\end{itemize}
It follows that the value of the two-pile game is $\max(k/n, 1-k/n)$, and that a blind guessing is the minimax strategy in asymmetric game.
\footnote{
A simpler function $g_n(x)=(\delta/2n) \,x^{-n\mp \delta}$ could have  been used (with the same breakpoint as in (\ref{gnon})), but this leads to $h_n(a)=\frac{\delta}{2}(1\pm\frac{\delta}{n})a^{\mp\delta-1}$ depending on $n$,
hence losing universality of the above solution.}

To practically implement the strategy defined by (\ref{diamond}),
(\ref{gnon})  Bob needs to simulate  a sample value $a$ (possibly   with the aid of techniques from
\cite{Devroye}), use a standard uniform random numbers generator  for $n$ values $u_i$, and set $z_i=au_i$.
To play  the best-response strategy in the Bayesian game Alice may  calculate
the  posterior distribution of the unknown scale parameter using the observed  $x=z_1\vee\ldots\vee z_k$, but this 
update  cannot help improving upon the blind guessing in the asymmetric game, and will give only a minor advantage in the symmetric one.

\end{document}